# Multi-objective Optimal Reactive Power Dispatch of Power Systems by Combining Classification Based Multi-objective Evolutionary Algorithm and Integrated Decision Making

**Meng Zhang, Yang Li, Senior Member, IEEE**

School of Electrical Engineering, Northeast Electric Power University, Jilin 132012, China

Corresponding author: Yang Li (e-mail: liyang@neepu.edu.cn).

This work was supported in part by the China Scholarship Council (CSC) under Grant 201608220144, the "13th Five-Year" Scientific Research Planning Project of Jilin Province Department of Education under Grant No. JJKH20200113KJ, and the National Natural Science Foundation of China under Grant No. 51677023.

**ABSTRACT** For the purpose of addressing the multi-objective optimal reactive power dispatch (MORPD) problem, a two-step approach is proposed in this paper. First of all, to ensure the economy and security of the power system, the MORPD model aiming to minimize active power loss and voltage deviation is formulated. And then the two-step approach integrating decision-making into optimization is proposed to solve the model. Specifically speaking, the first step aims to seek the Pareto optimal solutions (POSs) with good distribution by using a multi-objective optimization (MOO) algorithm named classification and Pareto domination based multi-objective evolutionary algorithm (CPSMOEA). Furthermore, the reference Pareto-optimal front is generated to validate the Pareto front obtained using CPSMOEA; in the second step, integrated decision-making by combining fuzzy c-means algorithm (FCM) with grey relation projection method (GRP) aims to extract the best compromise solutions which reflect the preferences of decision-makers from the POSs. Based on the test results on the IEEE 30-bus and IEEE 118-bus test systems, it is demonstrated that the proposed approach not only manages to address the MORPD issue but also outperforms other commonly-used MOO algorithms including multi-objective particle swarm optimization (MOPSO), preference-inspired coevolutionary algorithm (PICEAg) and the third evolution step of generalized differential evolution (GDE3).

## I. INTRODUCTION

Optimal reactive power dispatch (ORPD) of power systems refers to adjusting the parameters of the control equipment in the system to make the whole network at the optimal operation, which is of great significance to the economic and secure operation of the power system [1]. In recent years, with the improvement of the operation level of the power system, ORPD has evolved from a single-objective optimization problem to a multi-objective optimization (MOO) problem that comprehensively considers various operation indicators [2].

So far, there have been a large number of studies carried out to solve MOO issues. In some literature, the multi-objective optimal reactive power dispatch (MORPD) problem is transformed into a single-objective optimization problem by the scalarization reflecting the preference degree of each objective in advance. The arguably common scalarization approach termed weighted sum technique is adopted in [3, 4]. But the weight factors reflecting the decision-makers' preference are often difficult to determine in reality. Moreover, the traditional optimization approach can find one solution at most in a single simulation run, which makes the computation cost heavier [5]. Taking this into account, multi-objective evolutionary algorithms (MOEAs) have been used to solve the MORPD model [6]. Although the Pareto optimal solutions (POSs) can be obtained, it is difficult to determine the best compromise solutions since different decision-makers have different preferences for a given operational condition. Furthermore, the preference of the same decision marker can vary with the changing operating requirements of the system.

As we all know the objective functions are always conflict and can't be optimal at the same time, however, POSs which make a compromise of conflict objectives can be obtained after optimization. In [7] although the ORPD problem is treated as the MOO problem, there is no decision analysis regarding how to extract the best compromise solutions (BCSs) from the obtained POSs. Different from [7] decision-making step is considered in [5], however, there is only one BCS determined by the min-max criterion, which can't reflect the preference of decision-makers. In terms of the POSs, how to determine BCSs reflecting decision-makers' preferences is a significant and challenging problem. To address this issue, in [8] a two-stage approach taken the preference of decision-makers into consideration is put forward to overcome the issue in combined heat and power economic emission dispatch. In [9] the same approach is applied to handle hybrid AC/DC grids with voltage source converter based high voltage direct current problem. A two-stage optimization approach incorporating multi-objective optimization and decision analysis was employed to deal with distributed generation planning issues in distribution networks in [10]. In [11] a two-step approach is proposed to address a practical multi-objective dynamic optimal dispatch model for isolated micro-grids. In this paper, the two-step approach combining MOEAs with integrated decision-making analysis is proposed in this paper, where MOEAs have the ability to find the POSs in one run and decision-making analysis can determine the BCSs reflecting decision-makers' preferences from the POSs.

As for the MORPD issue, from the perspective of system security and economy, it is hoped that the active power loss will be minimized to reduce the investment and the voltage stability will be maximized to guaranteed voltage quality [12]. ORPD is a large-scale nonlinear mixed integer programming problem with continuous and discrete variables while satisfying both equality and inequality constraints [13]. Due to the powerful ability of MOEAs to find widely optimal POSs by using only one simulation run, MOEAs are widely used to solve the MORPD model. While MOEAs include enormous algorithms such as water cycle algorithm (NGBWC) [14], backtracking search optimizer (BSO) [15], whale optimization algorithm (WOA) [16], and grey wolf optimizer (GWO) [17], etc. Furthermore, for better performance, some improvements have been made based on the original algorithms. For example, in [18] a modified differential evolution algorithm (MDEA) is put forward to solve the ORPD problem to decrease the active power loss and voltage deviation. Similar to the approach in [18], a fuzzy adaptive heterogeneous comprehensive learning particle swarm optimization (PSO) algorithm is presented to address the MORPD problem through enhancing exploration and exploitation processes in [19]. For the same model, in [20] the improved gravitational search algorithm GSA-CSS based on conditional selection strategies (IGSA-CSS) improves the global search ability by using the memory characteristics of PSO, considering the shortcomings of GSA itself. In [21] an improved social spider optimization (ISSO) was proposed, through improvement the strong search ability was obtained due to the less value for control parameters. It can be seen that generally, the improvements involve how to enhance the global search ability, however, the improvements about this are rarely reported. In MOEAs, it is arguable that selection in evolutionary algorithms (EA) is essentially a classification problem on account of selection operators mainly based on objective values [22]. Following this idea, the classification based pre-selection (CPS) strategy was introduced into Pareto domination based MOEAs [22]. And further improvements were made in [23]. In this paper, a MOO algorithm named CPS base MEA (CPSMOEA) is introduced to solve the MORPD problem.

The main contributions of this work include the following aspects:

(1) To coordinate the economy and security of power systems, the CPSMOEA is introduced for the first time to solve the MORPD issue in this study.

(2) To determine the best comprise solutions from the Pareto optimal solutions, integrated decision-making analysis combining the fuzzy c-means algorithm (FCM) and grey relation projection method (GRP) is successfully employed.

(3) The simulation results on the IEEE 30-bus test system shown that the performance of the CPSMOEA is superior to that of multi-objective particle swarm optimization (MOPSO) and preference-inspired coevolutionary algorithm (PICEAg) in terms of convergence and distribution.

The rest of this article is structured as follows: the MORPD model is formulated in Section II. In section III, the model is solved by the proposed two-step approach incorporating the combination of FCM and GRP into the multi-objective optimization procedure. The simulation results on the IEEE 30-bus test system are given in section IV. Finally, the conclusion is drawn in section V.

## II. PROBLEM FORMULATION

As a typical MOO problem, the MORPD is formulated to achieve the ideal settings of control variables to satisfy certain objective functions, which can be described as:

$$\text{minimize } F(x,u) \qquad (1)$$

$$\text{subject to } \begin{cases} G(x,u) = 0 \\ H(x,u) \le 0 \end{cases} \qquad (2)$$

where $F(x,u)$ represents the objective function; $G(x,u)$ and $H(x,u)$ are respectively the equality and inequality constraints; $x$ and $u$ are respectively the vectors of dependent variables and control variables. In this study, the dependent variables refer to load bus voltages, while the control variables consist of generator bus voltages, transformer tap ratios and compensation capacity of shunt capacitor banks.

### A. OBJECTIVE FUNCTIONS
In this work, there are two objective functions: the active power loss and voltage deviation.

#### 1) ACTIVE POWER LOSS
For the economic view, with the reform of the power market, power suppliers always make the best use of the existing transmission capacity and active power. Therefore, reducing the active power loss on transmission lines has become an important issue concerned by the power department. The transmission loss is regarded as an objective function, as follows:

$$F_1 = \min P_{loss} = \sum_{k=1}^{N_L} g_k \left[ V_i^2 + V_j^2 - 2V_i V_j \cos(\delta_i - \delta_j) \right] \qquad (3)$$

where $N_L$ is the number of transmission lines; $V_i$ and $V_j$ are the voltage magnitudes at bus $i$ and $j$; $\delta_i$ and $\delta_j$ are the voltage angles at bus $i$ and $j$, respectively; and $g_k$ is transfer conductance between bus $i$ and $j$

#### 2) VOLTAGE DEVIATIONS
Considering the secure operation of the modern power system, voltage instability has become a critical issue that must be confronted. The objective function expressed as voltage deviation is used to evaluate voltage instability of the power system by minimizing the sum of voltage deviation at each load bus. It is defined as follows:

$$F_2 = \min VD = \sum_{k=1}^{N_{load}} \left| \frac{V_k - V_k^{ref}}{V^{upper} - V^{lower}} \right| \qquad (4)$$

where $N_{Load}$ is the number of load buses; $V_k^{ref}$ is the reference voltage at the $k$ th load bus which can be set to 1.0 p.u.; $V^{upper}$ is the upper limit of load bus voltage; $V^{lower}$ is the lower limit of load bus voltage.

### B. EQUALITY CONSTRAINTS
For any operating condition of the power system, the following two equality constraints containing active power balance and reactive power balance should be met.

$$P_{Gi} = P_{Li} + U_i \sum_{j \in Ni} U_j (G_{ij} \cos \theta_{ij} + B_{ij} \sin \theta_{ij}) \qquad (5)$$

$$Q_{Gi} = Q_{Li} + U_i \sum_{j \in Ni} U_j (G_{ij} \sin \theta_{ij} - B_{ij} \cos \theta_{ij}) \qquad (6)$$

where $P_{Gi}$ and $Q_{Gi}$ are active and reactive power generation at the bus $i$ respectively; $P_{Li}$ and $Q_{Li}$ are the load active and reactive power at the bus $i$, respectively; $G_{ij}$ and $B_{ij}$ are the transfer conductance and susceptance between bus $i$ and bus $j$, respectively.

### C. FILE FORMATS FOR GRAPHICS
#### 1) GENERATOR CONSTRAINTS
The generator operating under any condition should be within its upper and lower limits. The minimum and maximum boundaries of voltage and reactive power output are given below:

$$V_{Gi}^{\min} \le V_{Gi} \le V_{Gi}^{\max} \qquad i = 1, 2, \cdots, N_G \qquad (7)$$

$$Q_{Gi}^{\min} \le Q_{Gi} \le Q_{Gi}^{\max} \qquad i = 1, 2, \cdots, N_G \qquad (8)$$

where $N_G$ is the number of generators.

#### 2) TRANSFORMER CONSTRAINTS
Transformer tap settings vary between the maximum value and minimum value as follow:

$$T_i^{\min} \le T_i \le T_i^{\max} \qquad i = 1, 2, \cdots, N_T \qquad (9)$$

where $N_T$ is the number of transformers.

#### 3) SHUNT CAPACITOR BANK CONSTRAINTS
The compensation capacity of shunt capacitor banks should be restricted to the upper and lower boundaries as below:

$$Q_{Ci}^{\min} \le Q_{Ci} \le Q_{Ci}^{\max} \qquad i = 1, 2, \cdots, N_C \qquad (10)$$

where $N_C$ is the number of shunt capacitor banks.

#### 4) LOAD VOLTAGE CONSTRAINTS
The load bus voltages should be maintained in a reasonable range as follow:

$$V_{Li}^{\min} \le V_{Li} \le V_{Li}^{\max} \qquad i = 1, 2, \cdots, N_{load} \qquad (11)$$

where $N_{load}$ is the number of loads.

#### 5) SECURITY CONSTRAINTS
The apparent power flow on every transmission line should be limited to its allowable range to avoid overload, which is given below:

$$S_{li} \le S_{li}^{\max} \qquad i = 1, 2, \cdots, N_{line} \qquad (12)$$

where $N_{line}$ is the number of the transmission lines.

### III. MODEL SOLUTION

### A. SOLUTION FRAME

This section will primarily discuss the solution step of the MORPD model. When solving the single-objective optimization problem, a unique optimal solution is got. However, POSs can be obtained to balance conflict objectives in solving the MOO problem. Furthermore, for the decision-makers' preference, a two-step approach that contains the decision analysis is proposed. This approach is mainly divided into two steps: one is the optimization step and the other is the decision-making step.

Optimization step: A set of POSs can be obtained by solving the MORPD model formulated in section II by the MOO algorithm, i.e., CPSMOEA. Moreover, the reference Pareto-optimal front is obtained through multiple runs of single objective optimization with the weighted sum of objectives to validate the PF generated by the CPSMOEA. Here, particle swarm optimization approach based on multiagent systems (MAPSO) is adopted for performing single-objective optimization.

Decision-making step: First in this step, the POSs are divided into two clusters which represent the two different preferences of decision-makers via FCM, and then BCSs are respectively extracted from the two clusters via GRP.

The solution framework of the two-step approach is shown in Fig.1.

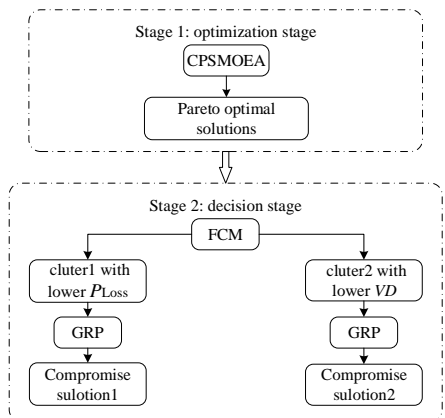

**FIGURE 1.** Solution framework.

## B. MULTI-OBJECTIVE OPTIMIZATION

### 1) MOEA

In EAs, the population is evaluated by the objective fitness function that indicates their pros and cons [24]. For the single-objective optimization problem, a unique objective function is used to evaluate the pros and cons, however, in the MOO problem the evaluation of the individual's pros and cons becomes a difficult issue due to the different objectives existed simultaneously. Goldberg proposed a new idea combining Pareto theory with EAs to solve the MOO problem, which is of great significance for the subsequent research on multi-objective evolutionary algorithms [25]. And further MOEAs gradually developed an algorithm with good practicability and robustness, which has received extensive attention [26]. It can be said that the introduction of the Pareto theory is the key to the wide application of MOEAs [27, 28].

Furthermore, the POSs consist of a set of solutions illustrating the regularity in both objective and decision spaces. Basing on Pareto domination, CPS is readily combined into Pareto domination based on MOEAs. Motivated by this idea, CPSMOEA is proposed to solve the MORPD problem.

### 2) CPS STRATEGY

In fact, CPS is essentially a classification problem. The purpose of the classification study is to predict the class labels of those unknown instances based on some known training data, mainly to extract as much information as possible from the known data. Unlike these existing methods which take classification as surrogate procedures, classification was employed to pre-selection in this paper.

In terms of pre-selection, there are many different meanings in MOEA [29]. In this paper, pre-selection refers to a procedure in which the promising solutions are extracted after the current solution generates candidate offspring solutions by reproduction operators while the rest unpromising ones are deserted. After the pre-selection procedure, the promising ones which selected as the offspring population are chosen into the next generation after environment selection. In this paper, the specific procedure to implement CPS is mainly divided into three parts: the labeling data set, the classifier training, and the selection of promising offspring solutions. The specific processes of these parts are given in detail as follows:

(1) Labeled Population

Here, the current populations are first used as the training data sets and separated them into two classes. One of the two classes is in the external population P+ with label +1 denoting the 'promising' training dates. And the other is in the external population P− with label −1 representing the 'unpromising' training dates. Since the MOEA adopted in this paper is based on Pareto domination, the sorting scheme of Pareto domination can be naturally incorporated into the data classification. In other words, if a solution is non-dominated, the label of it is +1 ; otherwise, it's label is −1 .

The expression $Q = NDS(P, N)$ is denoted sorting scheme of Pareto domination, it means that the best $N$ solutions are stored in $Q$ finally [30]. First, the population $P$ is sorted into several parts that be ranked based on non-domination, among which, the population with the lower (better) rank is preferred. In each part, the population is non-dominated with each other. And then the individuals in each part are ordered by calculating the crowding distance, where the solution located in a lesser crowded region is preferred. For each individual, there are two attributes termed non-dmination rank ($i_{rank}$) and crowding distance ($i_{distance}$), respectively. They are given as follows:

$$i \prec_n j, if (i_{rank} < j_{rank})$$
$$or((i_{rank} = j_{rank}) \qquad (13)$$
$$and \ (i_{distance} > j_{distance}))$$

where $\prec_n$ called crowded-comparison represent the selection process.

According to the fitness values expressed as $\{F(x_1), \cdots, F(x_i), \cdots, F(x_N)\}$, the population expressed as $P = \{x_1, \cdots, x_i, \cdots, x_N\}$ is separated into two classes, i.e., 'promising' populations which stored in $P+$ and 'unpromising' populations which stored in $P-$.

(2) Classifier Training

The classifier based on the sorted training data, respectively, expressed as $P+ = \{x_1, \cdots, x_i, \cdots, x_{\frac{N}{2}}\}$ and $P- = \{x_1, \cdots, x_i, \cdots, x_{\frac{N}{2}}\}$ is established in this procedure. Among which $x_i$ represents an individual with a multi-dimensional control vector. For each individual $x_i$, there is a label $l \in \{+1, -1\}$ corresponding to it. The training classifier aims to establish a connection between an individual $x_i$ and the label of it.

In this paper, K-Nearest Neighbor (KNN) is used to determine the label of the individual $x_i$, which is given by

$$l = \begin{cases} +1, & if \ \sum_{i=1}^{K} C(x_i) \geq 0 \\ -1, & otherwise \end{cases} \qquad (14)$$

where $K$ represents the number of neighbors taken into account in determining the class, $x_i$ denotes the $i$th closest control vector, $C(x_i)$ is the real relationship of $x_i$.

(3) Offspring Selection

Offspring selection is essentially to choose the 'promising' offspring solutions with good quality from all offspring candidates by evaluating the fitness value. Furthermore, how to produce candidate offspring solutions is an important component of this procedure. Herein, differential evolution (DE) reproduction operator is used to solve this problem. DE algorithm is an efficient population-based heuristic algorithm. DE reproduction operator generates candidate offspring solutions by mutation, crossover, and selection operation. The three operations are described as follow:

Mutation operation: In terms of each individual, there is a vector $P = \{x_1, \cdots, x_i, \cdots, x_N\}$ named target vector and the mutation operation produces a corresponding vector $V = [v_1, \cdots v_i, \cdots, v_N]$ termed donor vector. The basic idea of the mutation operation base on the DE is to add a difference vector to the base vector. The original mutation operation is given by

$$v_i = x_{r_1} + F \times (x_{r_2} - x_{r_3}) \qquad (15)$$

where $v_i \in V$ ; $N$ is the number of the individuals of the population; $r_1$, $r_2$ and $r_3$ are not equal to $i$, and they are random different integers in the interval $[1 \ N]$ ; the mutation control parameter $F$ is a positive number and usually limited to the interval $[0 \ 1]$.

Crossover operation: Crossover operation generates offspring individuals by performing discrete recombination on target vector and donor vector, that is, $U = [u_1, \cdots u_i, \cdots, u_N]$ termed trial vector. The basic idea of crossover operation is that donor vector and target vector exchange elements with each other to improve the population diversity. The specific implementation of crossover operation is as shown in Eq. (16):

$$u_{i,j} = \begin{cases} v_{i,j} & if(rand(0,1) \leq Cr) or(j = sn) \\ x_{i,j} & otherwise \end{cases} \qquad (16)$$

where $u_{i,j}$ represents the $j$th dimension element value meeting $u_{i,j} \in u_i = [u_{i,1}, \cdots u_{i,j} \cdots u_{i,n}]$ ; $sn$ is a random integer satisfying $sn \in [1, 2, \cdots, n]$ ; the crossover control parameter $Cr$ is a positive number and usually limited to the interval $[0 \ 1]$.

Selection Operation: The selection operation determines the evolution direction of the whole population, and the greedy choice is applied in this process. In terms of the parent individual $x_i$, if the corresponding offspring individual $u_i$ is worse than it, then the $x_i$ is selected for the next generation. Otherwise, the $u_i$ is selected and the selection process is as shown in Eq. (17):

$$x_i^{t+1} = \begin{cases} u_i^t & if \ f(u_i^t) \leq x_i^t \\ x_i^t & otherwise \end{cases} \qquad (17)$$

where $x_i^{t+1}$ represent parent individual in the next generation.

In pre-selection, the qualities of these candidate solutions are evaluated by means of CPS, and the promising one will be selected for the real function evaluation.

## C. APPLICATION OF CPSMOEA IN MORPD PROBLEM

Regarding the application of the CPSMOEA in solving MORPD problem, the main solution processes are as follows:

*Step 1*: Initialization of the system. Enter the following initial variables: 1) system parameters such as the data of buses, branches, loads and generators; 2) algorithm parameters such as the population size, the number of objectives and variables, and so on; 3) the boundaries and steps of related variables.

*Step 2*: Initialize individual vectors. The position of the individual in search space corresponds to the control variables. The continuous control variables are generator bus voltage $V_G$; while the discrete control variables comprise the transformer tap ratios $T$ and the compensation capacity $Q_C$ of shunt capacitor banks. The dimension of each individual vector is determined by the number of control variables expressed as

$$\underbrace{\left\{ V_{G,1}, \cdots, V_{G,N_G} \middle| T_1, \cdots, T_{N_T} \middle| Q_{C,1}, \cdots, Q_{C,N_C} \right\}}_{N}.$$

*Step 3*: Calculation of objective functions. For each individual, calculate the $P_{loss}$ and $VD$ according to equations (1) and (2).

*Step 4*: Classification of pre-selected data. Pareto domination is used to label the current population, and then the current population is sorted according to the rank crowding distance. The promising individuals with label $+1$ were stored in $P+$ and the rest unpromising with label $-1$ were stored in $P-$.

*Step 5*: The model of the classifier. Based on the pre-selected data in $P+$ and $P-$, KNN was used to find the relationship between each individual and the related label.

*Step 6*: Generate offspring solutions. First, the DE reproduction operator is used to generate the candidate offspring solution; and then select those with the label $+1$ as the offspring solution by the classifier.

*Step 7*: Environmental selection. The objective functions of the offspring population are calculated and the individual with higher fitness value was selected for the next iteration while the external population was updated, simultaneously.

The flow chart of the CPSMOEA algorithm is shown in Fig. 2.

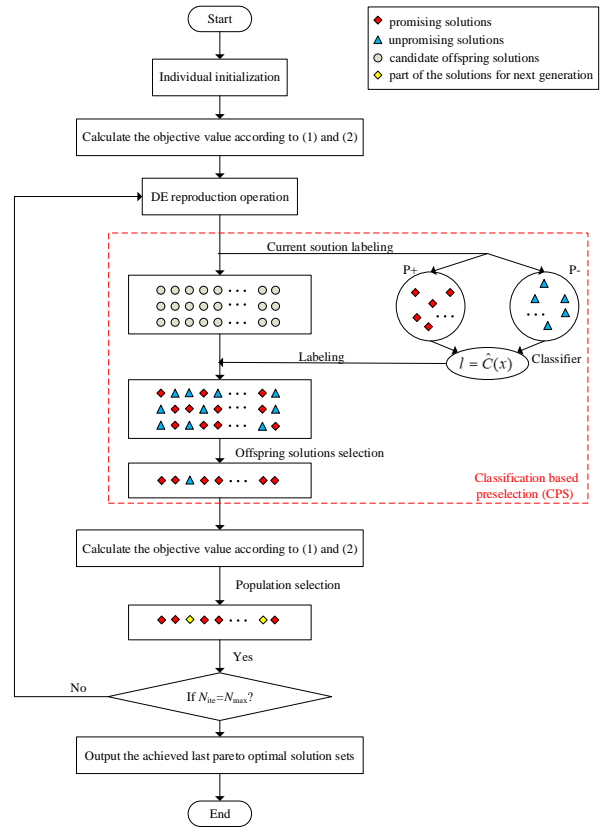

**FIGURE 2.** Flow chart of the CPSMOEA.

## D. DECISION SUPPORT

The Pareto optimal solution has a large scale, and the control vector contains different information. This paper proposes an auxiliary decision-making method combining FCM and GRP. It is convenient for the decision-makers to choose the compromise solution.

### 1) FUZZY C-MEANS CLUSTERING

FCM is an unsupervised clustering algorithm and its mathematical model is as follows:

$$\begin{cases} \min J_n(S,M,C) = \sum_{p=1}^{N_p} \sum_{q=1}^{N_{clu}} \eta_{p,q}^n \left\| s_p - c_q \right\|^2 \\ \mathrm{s.t.} \sum_{q=1}^{Nclu} \eta_{p,q} = 1 \end{cases} \quad (18)$$

where $J_n$ is a loss function; $S = \left\{ s_1, s_2, \cdots, s_p, \cdots, s_{N_p} \right\}$ is the vector of POSs; $N_p$ is the number of solutions in this vector; $M = \left\{ m_1, m_2, \cdots, m_{N_{clu}} \right\}$ and $C = \left\{ c_1, c_2, \cdots, c_q, \cdots, c_{N_{clu}} \right\}$ are respectively the membership degree matrix and cluster centers; $N_{clu}$ is a pre-given number of clusters; $\eta_{p,q}$ ($\eta_{p,q} \in [0,1]$) is membership degree representing the $S_p$ belongs to the $C_q$; $n$ ($n \in [1,\infty]$) is the fuzzy degree parameter.

Since two objectives are considered in this paper, the number of clusters is set to 2 for reflecting the different preferences of decision-makers over security and economy.

## 2) GREY CORRELATION PROJECTION

As an effective tool for handling various multiple attribute decision-making issues containing grey information, grey Correlation Projection (GRP) has been successfully used in many engineering areas [23, 24]. The projection of a scheme onto the ideal reference scheme is as follow:

$$pr_l^{+(-)} = \sum_{k=1}^{t} gr_{lk}^{+(-)} \frac{\lambda_k^2}{\sqrt{\sum_{k=1}^{t} \lambda_k^2}} \qquad (19)$$

where the superscript "$+$" indicates the ideal solution; and the superscript "$-$" indicates the negative ideal solution; $t$ is the total number of indicators; $gr_{lk}^{+(-)}$ is the grey correlation coefficient of the $k$ th indicator of the $l$ th scheme; $\lambda_k$ is the weight of each indicator of the scheme. In this paper, the weights corresponding to the two objectives are set to the same value, however, the decision-makers can adjust it according to the actual working condition or personal experience. The projection $p_l$ of each decision scheme on the ideal scheme is expressed as below:

$$p_l = \frac{pr_l^+}{pr_l^+ + pr_l^-} \qquad 0 \le p_l \le 1 \qquad (20)$$

where $p_l$ is the priority membership of scheme $l$ [8]. The scheme with the highest priority membership will be chosen as the BCSs.

### E. METRIC INDICATOR

For a variety of different algorithms, how to compare and measure their performance has become a meaningful topic. Many metric indicators are proposed to deal with the problem, the current mainstream metric indicators are mainly divided into three categories:

(1) Evaluate the degree of convergence of the POSs;

(2) Evaluate the distribution of the solutions over the whole Pareto front (PF), mainly considering uniformity and diversity.

(3) Evaluates the convergence and distribution of the solutions comprehensively.

Only using the pure index, the whole performers of the algorithm can't be reflected. However, the comprehensive index alone can't reflect the quality of the algorithm in a certain aspect, so this paper uses the above three metric indicators for comprehensive evaluation.

## 1) GENERATIONAL DISTANCE

Generational Distance (GD) is a credible metric indicator to evaluate the convergence of solutions by calculating the sum of solutions' adjacent distances [31]. It is defined as shown in Eq. (21).

$$GD = \frac{\sum_{x \in P^*} dist(x, S^*)}{|P^*|} \qquad (21)$$

where $x$ is an element in $P^*$ that represents the approximate solution set; $S^*$ is the set of targeted points; $dist(x, S^*)$ is the Euclidean distance between $x$ and the nearest individual that belongs to $S^*$. If the GD of the algorithm has a smaller value, it indicates that the algorithm has stronger convergence.

## 2) SPREAD

Spread is a pure metric indicator of distribution. In terms of distribution, it can assess both diversity and uniformity. The spread is defined as shown in Eq. (22).

$$spread = \frac{distD_f + distD_l + \sum_{i=1}^{N-1} \left| distD_i - \overline{distD} \right|}{distD_f + distD_l + (N-1)\overline{distD}} \qquad (22)$$

where $distD_f$ is the Euclidean distance of extreme solutions; $distD_l$ is the Euclidean distance of the boundary solutions [30]; $\overline{distD}$ is the average value of all distance $distD_i (i = 1, 2, \cdots N-1)$; $N$ is the number of the final non-dominated points. And the smaller the value of spread, the better the distribution of obtained PF [8].

## 3) INVERTED GENERATIONAL DISTANCE

Considering both convergence and distribution, inverted generational distance (IGD) is used to estimate the performance of the algorithm. IGD is widely used in the MOO problem due to its high computational efficiency. It is defined as shown in Eq. (23).

$$IGD = \frac{1}{S^*} \sum_{\eta \in S^*} \min \left\{ dist(\eta, P^*) \right\} \qquad (23)$$

where $\eta$ is an element in $S^*$ that represents the set of targeted points; $P^*$ is the approximate solution set; $dist(x, y)$ is the nearest distance from $\eta$ to $P^*$. It is noteworthy that the final obtained solution set with smaller values of IGD has better diversity and convergence [32].

## IV. NUMERICAL RESULTS

The proposed approach has been implemented to address the MORPD issue on IEEE 30-bus and IEEE 118-bus test systems.

### A. IEEE 30-BUS TEST SYSTEM

As a widely-used test system in the ORPD field, the well-known IEEE 30-bus test system is adopted in this paper to examine the effectiveness and superiority of the proposed approach [18-20]. The single-line diagram of this test system is shown in Fig. 3.

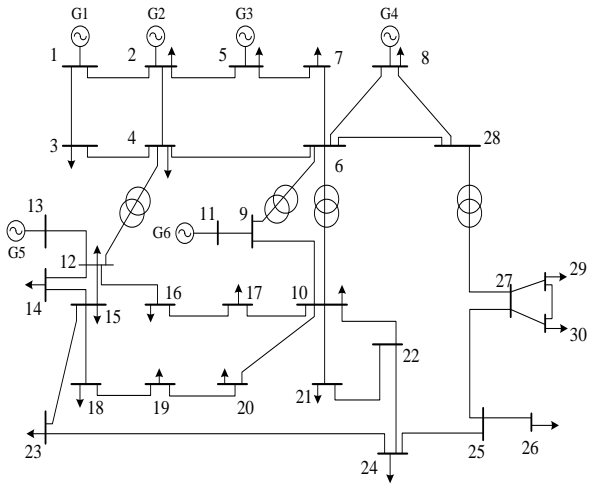

**FIGURE 3.** IEEE 30-bus test system.

As shown in Fig. 3, the system has 41 branches, 6 generators, and 22 loads [21]. The voltage amplitude of each generator is within the interval [0.9, 1.1], and the amplitude of each load bus voltage is within the interval [0.95, 1.05]. The four branches 6-9, 6-10, 4-12, and 27-28 are equipped with under-load tap-changing transformers which their taps vary in the range [0.9,1.1], and the step is 0.01 p.u. Shunt capacitor banks are installed on buses 3, 10 and 24. The number of shunt capacitor banks is 20, and the compensation capacity of each bank is 1 Mvar.

### 1) PARAMETER SETTINGS

In this paper, the control variables are encoded in a hybrid coding scheme. Specifically, the continuous control variables including generator bus voltage are real coded and the discrete variables comprised of the transformer tap ratio and compensation capacity of shunt capacitor banks are integer coded.

The relevant parameters of the CPSMOEA are given in Table I

TABLE I
PARAMENT SETTING

| Parameter | meanings | Values |
|---|---|---|
| $N$ | population size | 100 |
| $M$ | number of objectives | 2 |
| $D$ | number of variables | 13 |
| $eval$ | maximum number of evaluations | 10000 |
| $Cr$ | crossover control parameter | 1 |
| $F$ | mutation control parameter | 0.5 |

### 2) OPTIMIZATION RESULTS

In order to appropriately estimating the performance of the CPSMOEA, the PICEAg and MOPSO are introduced as the comparison algorithms. To further verify the effectiveness of the CPSMOEA , the reference Pareto-optimal front is generated by multiple runs of single objective optimization using the MAPSO. If the solutions obtained using an algorithm are close to the reference Pareto-optimal front, the algorithm is said to be good [6]. In this paper, the reference Pareto-optimal front consists of 100 non-dominated solutions obtained from 100 independent runs. As shown in Fig. 4, the reference Pareto-optimal front and the PFs obtained by CPSMOEA, PICEAg and MOPSO are respectively given.

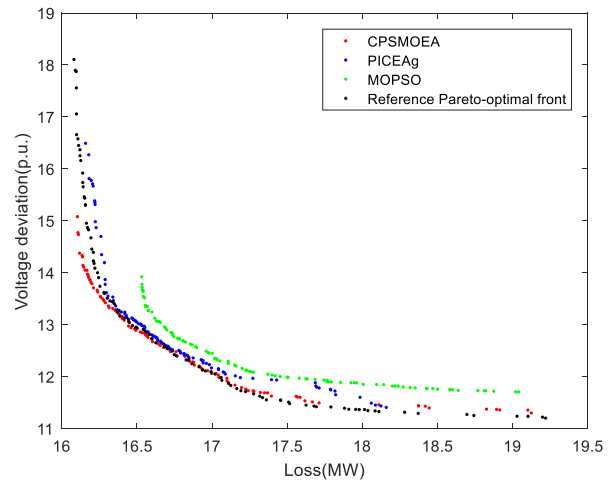

**FIGURE 4.** Reference Pareto-optimal front and the obtained PFs of CPSMOEA, PICEAg and MOPSO on IEEE 30-bus test system.

From Fig. 4, it can be easily found that the PF of the CPSMOEA dominates those of PICEAg and MOPSO, embodying that the curve obtained by the optimization of the CPSMOEA is closest to the two coordinate axes. Meanwhile, the PF obtained by the CPSMOEA is approaching the reference Pareto-optimal front more closely than that obtained by the PICEAg and MOPSO. This confirms the effectiveness and superiority of the presented method in terms of the optimization ability in the first step.

Moreover, it can be clearly seen that a more uniform PF can be obtained via the optimization of CPSMOEA. Hence a conclusion can be drawn that the optimization ability of the CPSMOEA in terms of distribution is better than those of the PICEAg and MOPSO.

To further examine the performance of the CPSMOEA, three different metric indicators GD, spread, and IGD are put forward. In view of the randomness of MOEAs [33], the proposed approach has been independently performed 30 times. After 30 independently simulation runs, the specific values of the three metric indicators are given in Table II.

TABLE II
THE VALUE OF METRIC INDICATORS GD, SPREAD, AND IGD

| algorithms | metrics | average value | best value | worst value |
|---|---|---|---|---|
| CPSMOEA | GD | 2.076 | 2.059 | 2.109 |
| | spread | 0.966 | 0.949 | 0.994 |
| | IGD | 20.458 | 20.437 | 20.479 |
| PICEAg | GD | 2.088 | 2.060 | 2.207 |
| | spread | 0.996 | 0.950 | 1.187 |
| | IGD | 20.501 | 20.464 | 20.599 |
| MOPSO [34] | GD | 2.129 | 2.071 | 2.249 |
| | spread | 0.978 | 0.963 | 0.998 |
| | IGD | 20.672 | 20.528 | 20.947 |
| GDE3 [35] | GD | 2.078 | 2.060 | 2.112 |
| | spread | 0.973 | 0.950 | 1.008 |
| | IGD | 20.462 | 20.467 | 20.485 |

From Table II, regarding the GD metric indicator, the average value of the CPSMOEA is 0.012 lower than that of the PICEAg, 0.053 lower than that of the MOPSO and 0.002 lower than that of the GDE3; the best value of the CPSMOEA is 0.001 lower than that of PICEAg, 0.012 lower than that of the MOPSO and 0.001 lower than that of the GDE3; the worst value of the CPSMOEA is 0.098 lower than that of the PICEAg, 0.140 lower than that of the MOPSO and 0.003 lower than that of the GDE3. GD is a metric indicator for evaluating convergence and the smaller the value of GD, the better the convergence of the algorithm. Consequently, it is demonstrated the optimization ability of the CPSMOEA in convergence performance is superior to that of the PICEAg, MOPSO, and GDE3.

In terms of the spread metric indicator, the average value of the CPSMOEA is equal to that of the PICEAg, 0.012 less than that of the MOPSO, 0.007 less than that of the GDE3; the best value of the CPSMOEA is 0.001 less than that of the PICEAg, 0.014 less than that of the MOPSO and 0.001 less than that of the GDE3; the worst value is 0.193 less than that of the PICEAg, 0.04 less than that of the MOPSO and 0.014 less than that of the GDE3. This proves that the distribution of PF obtained through the CPSMOEA is superior to those of PF obtained through three comparison algorithms.

Regarding the IGD metric indicator, compared with the PICEAg, MOPSO and GDE3 algorithms, the average value of the CPSMOEA is respectively decreased by 0.043, 0.214 and 0.004; the best value of the CPSMOEA is respectively decreased by 0.027, 0.091 and 0.03; the worst value of the CPSMOEA is respectively decreased by 0.12, 0.468 and 0.006. IGD is a metric indicator that considers convergence and distribution together. Thus, the above analysis demonstrates that in terms of convergence and distribution the CPSMOEA performs better than the PICEAg, MOPSO, and GDE3.

### 3) DECISION-MAKING ANALYSIS
In the first step of the proposed approach, the CPSMOEA is used to solve the formulated MORPD model. Fig. 5 shows the PF obtained by the CPSMOEA.

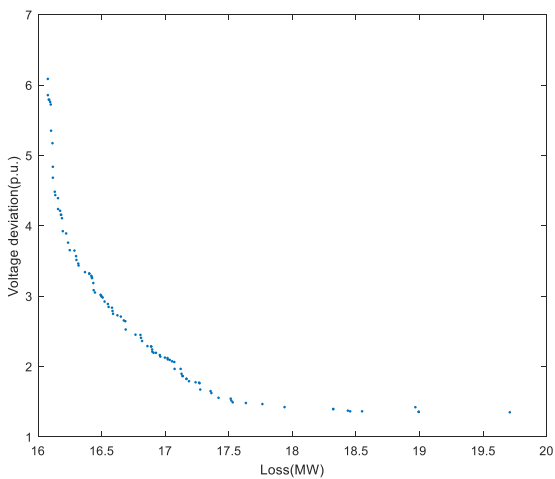

**FIGURE 5.** PF distribution of the CPSMOEA before clustering.

From Fig. 5 it can be clearly seen that as the active power loss decreases, the voltage deviation index increases simultaneously. Since the two objectives, i.e. active power loss and voltage deviation are conflicting, it cannot be optimal at the same time. However taking the different preferences of decision-makers' into account, it is not enough to obtain the POSs which simply reflecting the compromise. Thus, in the decision-making step of the proposed approach, integrated decision-making combining FCM and GRP is used to deal with this issue. After the second step, i.e., the decision-making step, the result is given in Fig. 6.

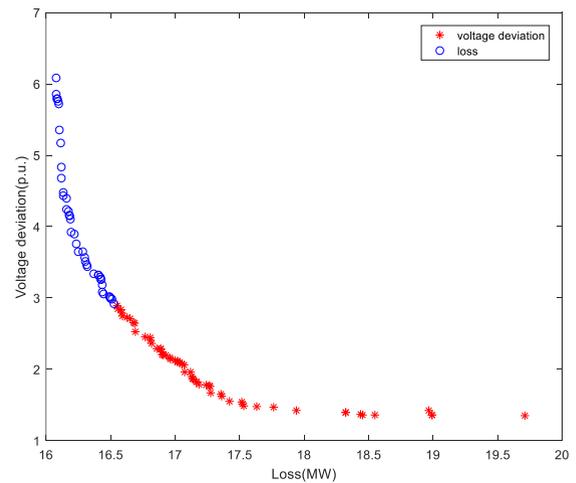

**FIGURE 6.** PF distribution of the CPSMOEA after clustering.

As shown in Fig. 6, the obtained POSs are divided into two clusters through the FMC, in this way, the different preferences of decision-makers are fully considered. And then, two BCSs are respectively chosen from each cluster via GRP in the decision-making step.

In order to properly evaluate the performance of integrated decision making, the BCSs obtained using the proposed method and the other comparison algorithms are shown in Table III. Note that here, the BCSs obtained by the other three algorithms only use the GRP (without consideration of FCM) in the decision making step.

TABLE III
BCSs OF DIFFERENT ALGORITHMS ON IEEE 30-BUS SYSTEM

| Items | Proposed method | | CPSM OEA | PICEAg | MOPSO [34] | GDE3 [35] |
|---|---|---|---|---|---|---|
| | BCS1 | BCS2 | BCS | BCS | BCS | BCS |
| $P_{loss}$ (MW) | 16.17 | 17.20 | 16.39 | 16.45 | 16.71 | 16.47 |
| $VD$ (p.u.) | 3.93 | 1.76 | 3.25 | 3.37 | 3.88 | 3.49 |

In Table III, BCS 1 and BCS 2 are obtained from two separate clusters after clustering, respectively. In terms of the preference on the economy, BCS 1 of the CPSMOEA is respectively decreased by 0.22 MV, 0.28 MV, 0.54 MV and 0.3 MV compare with the BCSs of CPSMOEA, PICEAg, MOPSO and GDE3; while its voltage deviation index is respectively increased by 0.68 p.u., 0.56 p.u., 0.05 p.u. and 0.44 p.u.. As far as the preference on the security is concerned, compared with the BCSs of the CPSMOEA, PICEAg, MOPSO and GDE3, BCS 2 of the proposed approach is respectively decreased by 1.49 p.u., 1.61 p.u. 2.12 p.u. and 1.73p.u.; while its power loss is increased by 0.81 MV, 0.75 MV, 0.49 MV and 0.73 MV. Furthermore, through the analysis of the specific data in TABLE III, it can be drawn that the BCSs obtained by GRP fully consider the preferences of decision-makers and provides more choices for decision-makers.

From the above analysis, the BCSs can be obtained through GRP. In the process of extracting BCSs, the GRP is first used to calculate the priority membership of the POSs, and then the solution with the highest priority membership is selected as BCS. In order to reasonably evaluate the optimization performances of the presented approach, a comparison test between the CPSMOEA and seven other algorithms, i.e., VaEA, NSGAIII, tDEA, IBEA, BiGE, MOPSO, and KnEA, has been performed. Here, the priority memberships of these algorithms are calculated by using the GRP. The comparison results are shown in Fig. 7.

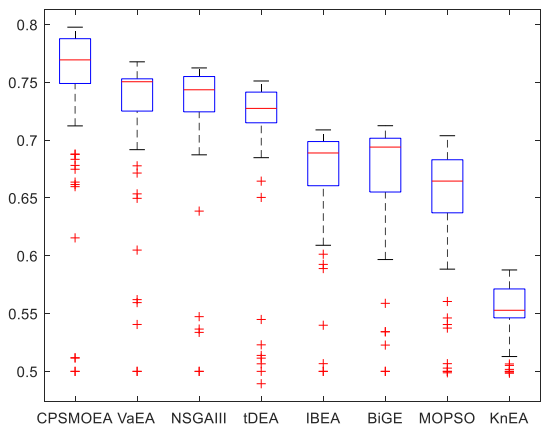

**FIGURE 7. Box plot of priority membership.**

Fig.7 shows the data distribution of the priority membership in a box plot. As mentioned in section Ⅲ, a greater value of priority membership signifies a better scheme will be obtained. Regarding the median of the priority membership, the value of CPSMOEA is maximum and the overall data distribution is superior to the other seven algorithms except the extreme solutions. Hence, a conclusion may be safely drawn based on the evidence that high-quality solutions can obtain via CPSMOEA to provide decision-makers with better choice of the BCSs point.

### 4) DISCUSSIONS OF BEST COMPARISON SOLUTIONS
To verify the availability of the proposed two-step approach, three different cases are taken into account. Before optimization, a unique solution can acquire corresponding to active power loss and voltage deviation on the basis of the original variables. After an independent simulation run, 100 Pareto optimal solutions and the corresponding individuals can be obtained, among which the comparison results of three different representative solutions are shown in Table IV.

Initial solution -- before optimization, a solution can be obtained.

BCS 1 -- after optimization the best compromise solution 1 reflects the decision-makers' preferences on the economy.

BCS 2 -- after optimization the best compromise solution 2 reflects the preferences of decision-makers on the security.

Through the two-step approach, the three different solutions are employed as the reference solutions, and the equivalent control variables are listed in Table IV.

**TABLE IV**
**COMPARATIVE RESULTS OF THREE DIFFERENT SOLUTIONS**

| Variables | Before optimization | After optimization | |
|---|---|---|---|
| | Initial solution | BCS 1 | BCS 2 |
| $V_{G1}$ (v) | 1.0600 | 1.1000 | 1.0836 |
| $V_{G2}$ (v) | 1.0430 | 1.0778 | 1.0530 |
| $V_{G5}$ (v) | 1.010 | 1.0417 | 1.0070 |
| $V_{G8}$ (v) | 1.0100 | 1.0478 | 1.0065 |
| $V_{G11}$ (v) | 1.0820 | 1.0393 | 0.9923 |
| $V_{G13}$ (v) | 1.0710 | 1.0293 | 1.0234 |
| $T_{6-9}$ (p.u.) | 0.98 | 1.05 | 1.01 |
| $T_{6-10}$ (p.u.) | 0.97 | 1.05 | 0.95 |
| $T_{4-12}$ (p.u.) | 0.93 | 1.05 | 0.98 |
| $T_{27-28}$ (p.u.) | 0.97 | 1 | 0.96 |
| $Q_{C3}$ (Mvar) | 5 | 12 | 1 |
| $Q_{C10}$ (Mvar) | 19 | 20 | 16 |
| $Q_{C24}$ (Mvar) | 4 | 12 | 14 |
| $P_{loss}$ (MW) | **17.46** | **16.17** | 17.20 |
| $VD$ (p.u.) | **6.38** | 3.93 | **1.76** |

Table IV shows that the proposed two-step approach can offer a reasonable solution scheme to the decision-makers according to their different preferences. Specifically speaking, $P_{loss}$ in BCS 1 through optimization is 7.4% lower than that in initial solution before optimization, meanwhile, $P_{loss}$ in BCS 2 after optimization is 1.5% lower than that in initial solution before optimization. The cause of this result is that the $P_{loss}$ of the two BCSs is both reduced through the optimization of the proposed approach while the reduction of $P_{loss}$ in BCS1 is greater due to the preferences on the economy. Similarly, $VD$ in BCS 2 after optimization is far less than that in initial solution before optimization (72.4%) and $VD$ in BCS 1 after optimization is less than that in initial solution before optimization (38.4%). It can be clearly seen from the above data analysis that the security of the power system has been improved after optimization. Furthermore, the value of $VD$ in BCS 2 has fallen still further, as BCS 2 pays particular attention to the security of the power system. As a result, a conclusion can be made that the distribution of power flow has a more reasonable trend.

### 5) COMPUTATIONAL EFFICIENCY ANALYSIS
To properly evaluate the efficiency of the CPSMOEA, the comparison tests between the CPSMOEA and three other algorithms PICEAg, MOPSO and GDE3 have been carried out. Accordingly, the test results are showed in Table V. Note that, considering the randomness of intelligent optimization algorithms, the average computational time of each algorithm in 30 independent runs is used as the computational time in the table.

TABLE V
COMPUTATIONAL TIMES OF THE ALGORITHMS

| Algorithms | Computational time (s) |
|---|---|
| **CPSMOEA** | **138.38** |
| PICEAg | 143.27 |
| MOPSO[34] | 153.65 |
| GDE3[35] | 143.81 |

As can be seen in Table V, the computational efficiency of the CPSMOEA is superior to that of the others. Moreover, the calculating efficiency of our approach can be further improved by using more advanced computer hardware and optimized code. As a result, ones can see that, in terms of computational efficiency, our method is able to meet the real-time requirements in practical applications, and that it precedes the other alternatives used in this study.

### B. IEEE 118-BUS TEST SYSTEM

To further establish the superiority of the proposed approach, IEEE 118-bus test system is adopted for an in-depth investigation. IEEE 118-bus test system consists of 54 generators, 9 transformers, 186 branches and 15 shunt capacitor banks [2, 21]. And the shunt capacitor banks are installed on 15 buses, respectively. Hence, the total number of control variables is 78. The upper and lower limits of the control variables on IEEE 118-bus test system are the same as those on IEEE 30-bus test system.

#### 1) OPTIMIZATION RESULTS

After optimization, the reference Pareto-optimal front, the results of the CPSMOEA and two comparison algorithms are plotted in Fig. 8, respectively.

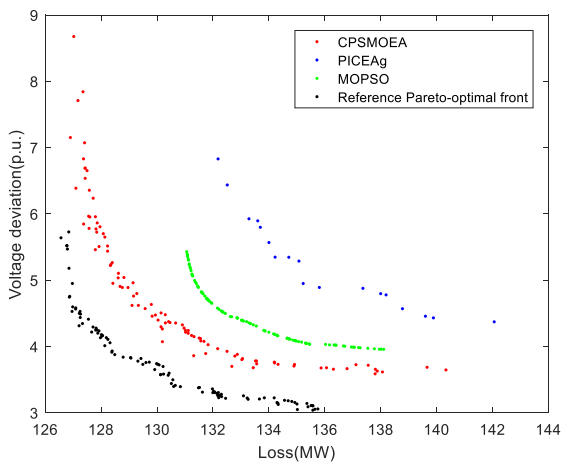

FIGURE 8. Reference Pareto-optimal Front and the obtained PFs of CPSMOEA, PICEAg and MOPSO on IEEE 118-bus test system.

From Fig. 8, it can be seen the PF obtained by the CPSMOEA is closer to the reference Pareto-optimal front than that obtained by the two comparison algorithms. And the PF obtained using the CPSMOEA dominates those obtained via the PICEAg and MOPSO. This shows that the CPSMOEA still outperforms the other two algorithms on IEEE 118-bus test system.

To make this result more convincing, three indicators consist of GD, spread and IGD are adopted for the comparison between the CPSMOEA and the other three algorithms. The average value, the best value and the worst value of the three metric indicators are detailed in Table VI.

TABLE VI
THE VALUE OF METRIC INDICATORS GD, SPREAD, AND IGD

| algorithms | metrics | average value | best value | worst value |
|---|---|---|---|---|
| CPSMOEA | GD | 13.169 | 12.886 | 13.458 |
| | spread | 0.983 | 0.963 | 0.990 |
| | IGD | 127.526 | 125.263 | 130.875 |
| PICEAg | GD | 24.333 | 22.120 | 39.103 |
| | spread | 0.987 | 0.965 | 1.002 |
| | IGD | 129.041 | 125.864 | 133.457 |
| MOPSO [34] | GD | 13.508 | 13.108 | 13.849 |
| | spread | 0.997 | 0.995 | 0.999 |
| | IGD | 133.309 | 129.403 | 136.401 |
| GDE3 [35] | GD | 13.516 | 13.264 | 13.973 |
| | spread | 0.983 | 0.964 | 0.991 |
| | IGD | 130.370 | 128.307 | 132.954 |

It can be clearly seen from Table VI, the average value, the best value and the worst value of all the metric indicators of the CPSMOEA are minimum. This suggests that the three metric indicators of the CPSMOEA are all better than those of the PICEAg, MOPSO and GDE3. Overall, this further reveals the superiority of the CPSMOEA in terms of convergence and distribution.

#### 2) DECISION-MAKING ANALYSIS

After the integrated decision-making analysis, two BCSs will eventually be obtained, as shown in TABLE VII.

TABLE VII
BCSs OF DIFFERENT ALGORITHMS ON IEEE 118-BUS SYSTEM

| Items | Proposed method | | CPSM OEA | PICEAg | MOPSO [34] | GDE3 [35] |
|---|---|---|---|---|---|---|
| | BCS 1 | BCS 2 | BCS | BCS | BCS | BCS |
| $P_{loss}$ ( MW) | 127.37 | 132.67 | 130.19 | 135.23 | 130.79 | 134.71 |
| VD (p.u.) | 5.85 | 3.70 | 4.07 | 4.95 | 4.62 | 4.22 |

From TABLE VII, it can be seen the two BCSs can be obtained after the integrated decision-making analysis. Moreover, the active power loss of BCS 1 and the voltage deviation of BCS 2 are superior to those of the BCS obtained by the CPSMOEA without integrated decision-making analysis. This fully takes the different preferences of decision-makers into account concerning economy and security. At the same time, it can also be seen from TABLE VII, without the integrated decision-making analysis, both active power loss and voltage deviation of BCS obtained by the CPSMOEA are less than those of the PICEAg, MOPSO, and GDE3, which further verifies the superiority of the CPSMOEA.

#### 3) COMPUTATIONAL EFFICIENCY ANALYSIS

After 30 independent runs, the average computational times of these algorithms are listed in TABLE VIII.

TABLE VIII
COMPUTATIONAL TIMES OF THE ALGORITHMS

| Algorithms | Computational time (s) |
|---|---|
| **CPSMOEA** | **271.43** |
| PICEAg | 288.84 |
| MOPSO [34] | 280.83 |
| GDE3 [35] | 298.87 |

According to the statistical results in TABLE VIII, the computational efficiency of the proposed CPSMOEA outperforms those of the PICEAg, MOPSO, and GDE3.

## V. CONCLUSION

This paper comprehensively considers the economics and security of the power system and establishes a model based on the active power loss and voltage deviation. A two-step approach containing a novel algorithm named CPSMOEA in its optimization step was first applied to the MORPD field. In the decision-making step of the proposed approach, the BCSs are selected from the POSs by means of the combination of the FCM and GRP, which fully considered the decision-makers' preference and can provide more schemes for the decision-maker. Based on the simulation results, it can be seen that the CPSMOEA can obtain a well-distributed Pareto front after introducing the CPS, and its convergence and distribution characteristics are better than the commonly-used multi-objective optimization algorithms such as PICEAg, MOPSO, and GDE3.

Our future work will focus on considering a dynamic index, like short-term voltage stability margin, as the objective function to handle dynamic security problems of the system. In addition, Another interesting topic is to extend this work to potential applications in multi-objective optimal operation of a microgrid/integrated energy system with uncertain renewable generations [36, 37].

## REFERENCES

[1] P.P. Biswas, P.N. Suganthan, R. Mallipeddi, "Optimal reactive power dispatch with uncertainties in load demand and renewable energy sources adopting scenario-based approach," Applied Soft Computing, vol. 75, pp. 616-632, Feb. 2019.

[2] M. Basu, "Quasi-oppositional differential evolution for optimal reactive power dispatch," International Journal of Electrical Power & Energy Systems, vol.78, pp. 29-40, Jun. 2016.

[3] H. Singh, L. Srivastava, "Modified differential evolution algorithm for multi-objective VAR management," International Journal of Electrical Power & Energy Systems, vol. 55, pp. 731-740, Feb. 2014.

[4] L. Srivastava, H. Singh, "Hybrid multi-swarm particle swarm optimisation based multi-objective reactive power dispatch," IET Generation, Transmission & Distribution, vol. 9, no. 8, pp. 727-739, May. 2015.

[5] Y. Li, P. Wang, H.B. Gooi, J. Ye, L. Wu. "Multi-objective optimal dispatch of microgrid under uncertainties via interval optimization," IEEE Transactions on Smart Grid, vol. 10, pp. 2046-2058, Mar. 2019.

[6] K. Deb, "Multi-objective optimization using evolutionary algorithms," Chichester, U.K.: Wiley, 2001.

[7] K. Nuaekaew, P. Artrit, N. Pholdee, and S. Bureerat, "Optimal reactive power dispatch problem using a two-archive multi-objective Grey wolf optimizer," Expert Systems with Applications, vol. 87, pp. 79–89, Nov. 2017.

[8] Y. Li, J. Wang, D. Zhao, G. Li, and C. Chen, "A two-stage approach for combined heat and power economic emission dispatch: Combining multi-objective optimization with integrated decision making," Energy, vol. 162, pp. 237–254, Nov. 2018.

[9] Y. Li, Y. Li, G. Li, D. Zhao, and C. Chen, "Two-stage multi-objective OPF for AC/DC grids with VSC-HVDC: Incorporating decisions analysis into optimization process," Energy, vol. 147, pp. 286–296, Mar. 2018.

[10] Y. Li, B. Feng, G. Li, J. Qi, D. Zhao, and Y. Mu, "Optimal distributed generation planning in active distribution networks considering integration of energy storage," Applied Energy, vol. 210, pp.1073-1081, Jan. 2018

[11] Y. Li, Z. Yang, D. Zhao, H. Lei, B. Cui, S. Li, "Incorporating energy storage and user experience in isolated microgrid dispatch using a multi-objective model," IET Renewable Power Generation, vol. 13, no. 6, pp. 973-981, Apr. 2019

[12] R. He, G.A. Taylor, Y.H. Song, "Multi-objective optimal reactive power flow including voltage security and demand profile classification," International Journal of Electrical Power & Energy Systems, vol. 30, no. 5, pp. 327-336, Jun. 2008.

[13] S.M. Mohseni-Bonab, A. Rabiee, B. Mohammadi-Ivatloo, "Voltage stability constrained multi-objective optimal reactive power dispatch under load and wind power uncertainties: A stochastic approach," Renewable Energy, vol. 85, pp. 598-609, Jan. 2016.

[14] A.A. Heidari, R.A. Abbaspour, A.R. Jordehi, "Gaussian bare-bones water cycle algorithm for optimal reactive power dispatch in electrical power systems," Applied Soft Computing, vol. 57, pp. 657-671, Aug. 2017.

[15] A.M. Shaheen, R.A. El-Sehiemy, S.M. Farrag, "Integrated strategies of backtracking search optimizer for solving reactive power dispatch problem," IEEE Systems Journal, vol. 12, no. 1, pp. 424-433, Jun. 2016.

[16] K.B.O. Medani, S. Sayah, A. Bekrar, "Whale optimization algorithm based optimal reactive power dispatch: A case study of the Algerian power system," Electric Power Systems Research, vol. 163, pp. 696-705, Oct. 2018.

[17] I.B.M. Taha, E.E. Elattar, "Optimal reactive power resources sizing for power system operations enhancement based on improved grey wolf optimizer," IET Generation, Transmission & Distribution, vol. 12, no. 14, pp. 3421-3434, Aug. 2018.

[18] W.S Sakr, R.A El-Sehiemy, A.M Azmy, "Adaptive differential evolution algorithm for efficient reactive power management," Applied Soft Computing, vol. 53, pp. 336-351, Apr. 2017.

[19] E. Naderi, H. Narimani, M. Fathi, "A novel fuzzy adaptive configuration of particle swarm optimization to solve large-scale optimal reactive power dispatch," Applied Soft Computing, vol. 53, pp. 441-456, Apr. 2017.

[20] G. Chen, L. Liu, Z. Zhang, "Optimal reactive power dispatch by improved GSA-based algorithm with the novel strategies to handle constraints," Applied Soft Computing, vol. 50, pp. 58-70, Jan. 2017.

[21] T.T Nguyen, D.N Vo, "Improved social spider optimization algorithm for optimal reactive power dispatch problem with different objectives," Neural Computing and Applications, pp. 1-32, Feb. 2019.

[22] J. Zhang, A. Zhou, G. Zhang, "A classification and Pareto domination based multiobjective evolutionary algorithm," in Proc. IEEE Congress on Evolutionary Computation (CEC), Sendai, Japan, 2015, pp. 2883-2890.

[23] J. Zhang, A. Zhou, K. Tang, and G. Zhang, "Preselection via classification: A case study on evolutionary multiobjective optimization," Information Sciences, vol. 465, pp. 388-403, Oct. 2018.

[24] A. Santiago, B. Dorronsoro, A.J. Nebro, "A novel multi-objective evolutionary algorithm with fuzzy logic based adaptive selection of operators: FAME," Information Sciences, vol. 471, pp. 233-251, Jan. 2019.

[25] G. Chen, J. Qian, Z. Zhang, and Z. Sun, "Applications of novel hybrid bat algorithm with constrained Pareto fuzzy dominant rule on multi-objective optimal power flow problems," IEEE Access, vol. 7, pp. 52060-52084, Apr. 2019.

[26] P.C. Roy, K. Deb, and M.M. Islam, "An efficient nondominated sorting algorithm for large number of fronts," IEEE Transactions on Cybernetics, vol. 49, no. 3, pp. 859-869, Mar. 2019.

[27] Z. He, J. Zhou, L. Mo, H. Qin, X. Xiao, B. Jia, and C. Wang, "Multiobjective reservoir operation optimization using improved multiobjective dynamic programming based on reference lines," IEEE Access, vol. 7, pp. 103473-103484, Jul. 2019.

[28] Z.A. Khan, A. Khalid, N. Javaid, "Exploiting nature-inspired-based artificial intelligence techniques for coordinated day-ahead scheduling to efficiently manage energy in smart grid," IEEE Access, vol. 7, pp. 140102-140125, Sep. 2019.

[29] S.W. Mahfoud, "Crowding and preselection revisited," Parallel Problem Solving From Nature, vol. 2, pp. 27-36, Apr. 1992.

[30] K. Deb, A. Pratap, S. Agarwal, "A fast and elitist multiobjective genetic algorithm: NSGA-II," IEEE Transactions on Evolutionary Computation, vol. 6, no. 2, pp. 182-1972002, Aug. 2002.

[31] Y. Liu, J. Wei, X. Li, "Generational distance indicator-based evolutionary algorithm with an improved niching method for many-objective optimization problems," IEEE Access, vol. 7, pp. 63881-63891, May. 2019.

[32] Y. Sun, G.G. Yen, Z. Yi, "IGD indicator-based evolutionary algorithm for many-objective optimization problems," IEEE Transactions on Evolutionary Computation, vol. 23, no. 2, pp. 173-187, Jan. 2018.

[33] T. Han, Y. Chen, J. Ma, "Surrogate modeling-based multi-objective dynamic VAR planning considering short-term voltage stability and transient stability," IEEE Transactions on Power Systems, vol. 33, no. 1, pp. 622-633, Apr. 2017.

[34] B. Zhou, Chan K W, Yu T, "Strength pareto multigroup search optimizer for multiobjective optimal reactive power dispatch," IEEE Transactions on Industrial Informatics, vol. 10, no. 2, pp. 1012-1022, May. 2017.

[35] S. Ramesh, S. Kannan, S. Baskar, "An improved generalized differential evolution algorithm for multi-objective reactive power dispatch," Engineering Optimization, vol. 44, no. 4, pp. 391-405, Apr. 2012.

[36] Y. Li, Z. Yang, G. Li, D. Zhao, and W. Tian, "Optimal scheduling of an isolated microgrid with battery storage considering load and renewable generation uncertainties," IEEE Transactions on Industrial Electronics, vol. 66, no. 2, pp. 1565-1575, Feb. 2018.

[37] Y. Li, C. Wang, G. Li, J. Wang, D. Zhao, and C. Chen, "Improving operational flexibility of integrated energy system with uncertain renewable generations considering thermal inertia of buildings," Energy Conversion and Management, vol. 207, pp. 112526, March 2020.